\documentclass[12pt]{amsart}
\usepackage{hyperref}
\usepackage[utf8]{inputenc}
\usepackage{amssymb, graphicx}
\usepackage{amsfonts}
\usepackage{amsmath}
\usepackage{latexsym}
\usepackage{amscd}
\usepackage{verbatim}
\usepackage{mathrsfs}
\usepackage{color}

\allowdisplaybreaks[4]
\newcommand{\mf}{\mathfrak}
\newcommand{\g}{\mf{g}}
\newcommand{\h}{\mf{h}}

\newcommand{\gl}{\mf{gl}}
\newcommand{\p}{\mf{p}}

\newcommand{\Z}{{\mathbb Z}}
\newcommand{\C}{{\mathbb C}}
\newcommand{\N}{{\mathbb N}}

\newcommand\pt[1]{\frac{\partial}{\partial t_{#1}}}

\allowdisplaybreaks[4]

\def\sl{\mathfrak{sl}}
\def\s{\mathfrak{s}}
\def\p{\mathfrak{p}}
\def\b{\mathfrak{b}}

\def\span{\mathrm{span}}

\newtheorem{theorem}{Theorem}[section]

\newtheorem{lemma}[theorem]{Lemma}
\newtheorem{corollary}[theorem]{Corollary}

\theoremstyle{remark}

\numberwithin{equation}{section}

\begin{document}
\title[]{Generalized Verma modules over $\sl(m+1)$ induced from simple highest weight modules}
\author[]{Yaohui Xue$^1$,\ Yan Wang$^{2*}$  \bigskip\\
{$^1$ S\MakeLowercase{chool of} M\MakeLowercase{athematics and} S\MakeLowercase{tatistics}, N\MakeLowercase{antong} U\MakeLowercase{niversity}, J\MakeLowercase{iangsu}, P. R.  C\MakeLowercase{hina}}\\
{$^2$ S\MakeLowercase{chool of} M\MakeLowercase{athematics}, T\MakeLowercase{ianjin} U\MakeLowercase{niversity}, T\MakeLowercase{ianjin}, P. R. C\MakeLowercase{hina}}\\
{$^*$ C\MakeLowercase{orresponding author: wangyan09@tju.edu.cn}}
}
\maketitle
\begin{abstract}
A class of generalized Verma modules over $\sl(m+1)$ are constructed from simple highest weight $\gl(m)$-modules. Furthermore, the simplicity criterion for these $\sl(m+1)$-modules are determined and an equivalence between generalized Verma modules and tensor modules are established.

\vspace{0.3cm}

\end{abstract}

\section{Introduction}
Let $\g$ be a finite-dimensional complex simple Lie algebra and $U(\g)$ be its universal enveloping algebra. Fix a Cartan subalgebra $\h$ of $\g$ and a triangular decomposition $\g=\g_-\oplus\h\oplus\g_+$ of $\g$. Denote by $\Delta$ the root system associated to $(\g, \h)$, by $\Delta_+\subset \Delta$ the set of positive roots, and by $\Pi\subset \Delta_+$ the set of simple roots. Choose a subset $I \subset \Pi$ and denote by $\Delta_I$ the root subsystem generated by $I$. Let $\p=\mathfrak{l}\oplus \mathfrak{u}$ be the standard parabolic subalgebra corresponding to $I$, where $\mathfrak{l}=\h\oplus\sum_{\alpha\in\Delta_I}\g_\alpha$ is a reductive subalgebra of $\g$ and $\mathfrak{u}=\sum_{\alpha\in\Delta_+\setminus\Delta_I}\g_\alpha$ is the nilradical of $\p$.

Let $V$ be a simple $\p$-module with $\mathfrak{u}V=0$. Then the generalized Verma module (GVM) $M_\p(V)$ is the induced module
$$M_\p(V):=\mathrm{Ind}_{\p}^{\g}(V)= U(\g)\otimes_{U(\p)}V.$$
If $\p$ happens to be the Borel subalgebra $\b$ and $V$ is one-dimensional, then we obtain the classical Verma module given in \cite{V}. The concept of generalized Verma modules for finite-dimensional complex simple Lie algebras was  initially introduced by Garland and Lepowsky in \cite{GL}. One of the most important results was shown in \cite{CF, DMP, Fe, Fu1} that any simple weight (with respect to a fixed Cartan subalgebra) module over a finite-dimensional complex simple Lie algebra is either cuspidal or a quotient of a certain generalized Verma module. The structure of generalized Verma modules is also related to the conformal field theory in theoretical physics, see for example \cite{Se}.

One of the fundamental problems in the study of generalized Verma modules is to establish a clear and practical simplicity criterion. This problem has been studied in many cases. The case of generalized Verma modules induced from finite-dimensional modules was studied in \cite{J,R}. For generalized Verma modules induced from infinite-dimensional modules, a lot of results were also obtained, see for example \cite{FM,KM1,KM2,MaS,Mc,MS}. Notably, \cite{MaS} established an equivalent condition for the simplicity of a generalized Verma module $M_{\p}(V)$ for an arbitrary simple $\p$-module $V$, under the restriction that the reductive part of $\p$ is of type $A$. This work introduced the so-called ``rough structure" of generalized Verma modules. The above results gave many theoretical simplicity criteria for generalized Verma modules, but some of them are very complicated in application.

For any $m\in\N$, let $\gl(m)$ be the general linear Lie algebra consisting of all $m\times m$ matrices, and $\sl(m)$ be the special linear Lie algebra consisting of all traceless $m\times m$ matrices. Here we focus on the case where $\g=\sl(m+1)$ and $\p$ is the sum of $\sl(m)$ and the standard Borel subalgebra. Under this setting, the reductive part $\mathfrak l$ of $\p$ is isomorphic to $\gl(m)$. If $V$ is an $U(\h)$-free $\mathfrak l$-module of rank $1$, then the generalized Verma module $M_{\p}(V)$ was studied in \cite{CLNZ}, where the simplicity criterion for $M_{\p}(V)$ was given. This pair $(\g,\p)$ is of special importance due to the existence of an equivalence (with respect to an isomorphism of $\g$) between $M_{\p}(V)$ and a tensor module over the Witt algebra $W_m$ restricted to $\g$, as detailed in Section 5. This equivalence was also observed in \cite{GLZ} for a special case. From the GVM point of view, \cite{GLZ} further characterized the simplicity criterion for $M_{\p}(V)$, where $V$ a simple cuspidal $\mathfrak l$-module of rank $1$.

In this paper, we consider the case where $V=L(\mu)$ is the simple highest weight module with the highest weight $\mu$ satisfying $\mu_1=\dots=\mu_{\bar i},\ \mu_{\bar i+1}=\dots=\mu_m$ for some $1\leq \bar i\leq m$. We give the simplicity criterion for $M_\p(V)$. This $\mathfrak l$-module $V$, when restricted to $\sl(m)$, is the unique simple quotient of the scalar generalized Verma module, which has been studied in \cite{BJ,BXX} using the Gelfand–Kirillov dimensions. From the GVM point of view, the $\sl(m+1)$-module studied in \cite{TZ}, where $V$ is one-dimensional, is just a special instance of this paper. Certainly, the work presented in this paper also provides guidance for future research on the case where $\mu$ is an arbitrary highest weight.

The paper is organized as follows. In Section 2, we introduce the necessary notations and preliminaries. In Section 3, we discuss some properties of the simple highest weight modules in some special cases of $\gl(m)$, which serve as foundational results for subsequent sections. In Section 4, we give the necessary and sufficient condition for the generalized Verma module to be not simple. Finally, in Section 5, we establish an equivalence between a class of generalized verma modules and tensor modules when restricted to $\sl(m+1)$.

\section{Preliminaries}
We denote by $\N, \Z_+,\Z,\C$ and $\C^*$ the sets of all positive integers, non-negative integers, integers, complex numbers and nonzero complex numbers, respectively. All algebras and vector spaces are assumed to be over $\C$. Let $m$ be a positive integer throughout this paper. For $\gamma=(\gamma_1,\ldots,\gamma_m)\in\C^m$, we set $|\gamma|=\sum\limits_{i=1}^m\gamma_i$.

Suppose $A$ and $B$ are associative algebras or Lie algebras, $\varphi:A\rightarrow B$ is an algebra homomorphism, and $V$ is a $B$-module. Define the $A$-module $V^\varphi$ as follows: $V^\varphi=V$ as vector space, and
$$a\cdot v=\varphi(a)v,\ \forall a\in A, v\in V.$$

\subsection{The general linear Lie algebra}
Recall that $\gl(m)$ has a basis $\{E_{i, j}\mid 1\leq i, j \leq m\}$ and the Lie bracket:
$$[E_{i, j} , E_{k,l}]=\delta_{j, k}E_{i,l}-\delta_{l,i}E_{k, j}$$
for $1\leq i, j, k, l \leq m$. Note that $\gl(m)$ is reductive with the decomposition of ideals $\gl(m)=\sl(m)\oplus \C I_m$, where $I_m=\sum\limits_{i=1}^m E_{i,i}$ is the identity matrix.

We have the triangular decomposition $\gl(m)=\gl(m)_-\oplus\h\oplus\gl(m)_+$, where
$$\gl(m)_+=\span\{E_{i,j}\mid 1\leq i<j\leq m\},$$
$$\gl(m)_-=\span\{E_{i,j}\mid 1\leq j<i\leq m\},$$
and $\h=\span\{E_{i,i}\ |1\leq i\leq m\}$ is standard Cartan subalgebra of $\gl(m)$. Denote the dual vector space of $\h$ by $\h^*$. For $\mu\in\h^*$, we identify $\mu$ with the vector $(\mu_1,\ldots, \mu_m)\in \C^m$, where $\mu_i=\mu(E_{i,i})$. Denote by $e_i$ the vector with $1$ as its $i$-th coordinate and $0$ as all other coordinates. The set $\Phi^+=\{e_i-e_j \mid 1 \leq i< j \leq m\}$ forms the set of positive roots of $\gl(m)$.

A $\gl(m)$-module $V$ is called a weight module if $V=\oplus_{\mu\in\C^m}V_\mu$, where
$$V_{\mu}=\{v\in V\ |\ E_{i,i}\cdot v=\mu_iv,\ 1\leq i\leq m\}$$
is called the weight space with weight $\mu$. We know that $\h\oplus \gl(m)_+$ is the Borel subalgebra of $\gl(m)$. For any $\mu\in\C^m$, let $\C v_\mu$ be the one-dimensional $\h\oplus\gl(m)_+$-module defined by $\gl(m)_+\cdot v_\mu=0$ and $E_{i,i}\cdot v_\mu=\mu_i v_\mu$ for $1\leq i\leq m$. Then we have the Verma module
$$M(\mu)=U(\gl(m))\otimes_{U(\h\oplus \gl(m)_+)}\C v_\mu\cong U(\gl(m)_-)\otimes\C v_\mu.$$
Furthermore, $M(\mu)$ has a unique simple quotient module $L(\mu)$, which is the simple $\gl(m)$-module with highest weight $\mu$. In fact,
any finite-dimensional simple $\gl(m)$ module is isomorphic to $L(\mu)$ for some $\mu$ with $\mu_s-\mu_{s+1}\in\Z_+,\ 1\leq s \leq m-1$.

\subsection{The spcial linear Lie algebra}
It is well known that the set
$$\{E_{k,l}\mid 1\leq k\neq l\leq m+1\}\cup \{E_{i,i}-E_{i+1,i+1}\mid 1\leq i\leq m\}$$
is the standard basis of the Lie algebra $\sl(m+1)$, which is simply denoted by $\s$ in this paper. We identify any matrix in $\sl(m)$ with a matrix in $\sl(m+1)$ in the natural sense.

Let $\mathfrak l=\mbox{span}\{E_{i,j},  E_{i,i}-E_{m+1,m+1}\mid 1\leq i\neq j\leq m\}$. Obviously, $\mathfrak l$ is a subalgebra of $\sl(m+1)$. By letting
$$E_{i,j}\mapsto E_{i,j},\ E_{i,i}-E_{m+1,m+1}\mapsto E_{i,i}$$
for $1\leq i\neq j\leq m$, we have that $\mathfrak l$ is isomorphic to $\gl(m)$.
Let $\mathfrak u=\sum\limits_{s=1}^m\C E_{s,m+1}$. Then it is easy to see that $\p=\mathfrak l\oplus \mathfrak u$ is a maximal parabolic subalgebra of $\sl(m+1)$.

Let $L(\mu)$ be the simple highest weight module of $\gl(m)$ with the highest weight $\mu\in\C^m$. Set $E_{s,m+1}\cdot L(\mu)=0$ for $1\leq s\leq m$. Then $L(\mu)$ becomes a module over $\p$. Therefore, we can get the  generalized Verma module $M_{\p}(V)=\mathrm{Ind}_\p^{\sl(m+1)}L(\mu)$ over $\sl(m+1)$.

\subsection{Tensor modules over Witt algebra}\label{ssectensorm}
For any $\alpha=(\alpha_1,\dots,\alpha_m)\in\Z^m$, set $t^\alpha:=t_1^{\alpha_1}\cdots t_m^{\alpha_m}$. The Witt algebra $W_m$ is the derivation Lie algebra of the Laurent polynomial algebra $A_m=\C[t_1^{\pm 1},\dots,t_m^{\pm 1}]$ and has a basis
$\{t^\alpha\pt{i}\mid \alpha\in\Z^m,\ 1\leq i\leq m\}$ with the bracket defined by
$$[t^\alpha\pt{i},t^{\alpha'}\pt{j}]=t^\alpha\pt{i}(t^{\alpha'})\pt{j}-t^{\alpha'}\pt{j}(t^\alpha)\pt{i}.$$

The Weyl algebra $K_m$ is the unital associative algebra
$$\C[t_1^{\pm 1},\dots,t_m^{\pm 1},\pt{1},\dots,\pt{m}].$$
There is a Lie algebra homomorphism $\pi$ from $W_m$ to the tensor algebra $K_m\otimes U(\gl(m))$ given by
$$\pi(t^\alpha\pt{i})=t^\alpha\pt{i}\otimes 1+\sum_{s=1}^m\pt{s}(t^\alpha)\otimes E_{s,i},\ 1\leq i\leq m.$$
Let $P$ be a $K_m$-module and $W$ be a $\gl(m)$-module. Then we have the $W_m$-module $F(P,W):=(P\otimes W)^{\pi}$. This module is called a tensor module over $W_m$.

The tensor module serves as a natural generalization of the well-known Shen-Larsson module, which was independently constructed by Shen in \cite{S} and Larsson in \cite{La}. The simplicity of tensor modules over Witt algebras was fully characterized in \cite{LLZ}, and further related results can be found in the references cited therein.

\section{Simple highest weight $\gl(m)$-modules for some special cases}
In this section, we give some properties of the $\gl(m)$-module $L(\mu)$ for some special cases. From now on, $L(\mu)$ is written as $V$ and $\mu=(\mu_1,\ldots,\mu_m)$.
\begin{lemma}\label{Lem3.1}
If $\mu_i=\mu_{i+1}=\dots=\mu_j$ for $1\leq i<j\leq m$, then
$$V_{\mu-k_i(e_i-e_{i+1})-\dots-k_{j-1}(e_{j-1}-e_j)}=0$$
for any $k_i,\dots,k_{j-1}\in\Z_+$ with $k_i+\dots+k_{j-1}\neq 0$.
\end{lemma}
\begin{proof}
By the PBW theorem, any element in $V_{\mu-k_i(e_i-e_{i+1})-\dots-k_{j-1}(e_{j-1}-e_j)}$ is spanned by elements of the form $E_{i_1,j_1}\cdots E_{i_s,j_s}\cdot v_\mu$, where $i\leqslant j_l<i_l\leqslant j$ for any $1\leqslant l\leqslant s$ and $\sum_{l=1}^s (e_{j_l}-e_{i_l})=k_i(e_i-e_{i+1})+\dots+k_{j-1}(e_{j-1}-e_j)$.

Suppose $s=1$. We will show that $E_{i_1,j_1}\cdot v_\mu=0$ by an induction on $i_1-j_1$. Note that $E_{i_1,j_1}$ is a root vector of root which is a sum of $i_1-j_1$ negative simple roots and $E_{p,p+1}\cdot E_{i_1,j_1}\cdot v_\mu=[E_{p,p+1},E_{i_1,j_1}]\cdot v_\mu$ for any $p\in\{1,\dots,m-1\}$.

If $i_1-j_1=1$, then $E_{p,p+1}\cdot E_{i_1,j_1}\cdot v_\mu=\delta_{p,j_1}(\mu_{j_1}-\mu_{i_1})v_\mu=0$. So $\gl(m)_+\cdot E_{i_1,j_1}\cdot v_\mu=0$, which means that $E_{i_1,j_1}\cdot v_\mu$ is a highest weight vector. So $E_{i_1,j_1}\cdot v_\mu=0$.

If $i_1-j_1>1$ and $[E_{p,p+1},E_{i_1,j_1}]\neq 0$, then $[E_{p,p+1},E_{i_1,j_1}]$ is a root vector of root that is a sum of $i_1-j_1-1$ negative simple roots. By the induction hypothesis, $E_{p,p+1}\cdot E_{i_1,j_1}\cdot v_\mu=0$. So $\gl(m)_+\cdot E_{i_1,j_1}\cdot v_\mu=0$ and it follows that $E_{i_1,j_1}\cdot v_\mu=0$.

For $s>1$, $E_{i_1,j_1}\cdots E_{i_s,j_s}\cdot v_\mu=E_{i_1,j_1}\cdots E_{i_{s-1},j_{s-1}}(E_{i_s,j_s}\cdot v_\mu)=0$. So $E_{i_1,j_1}\cdots E_{i_s,j_s}\cdot v_\mu=0$ for any $s\in\N$. Thus, $V_{\mu-k_i(e_i-e_{i+1})-\dots-k_{j-1}(e_{j-1}-e_j)}=0$.
\end{proof}

\begin{lemma}\label{Lem3.2}
Suppose $\mu_{\bar i+1}=\dots=\mu_t$ for some $\bar i$ with $1\leq {\bar i}< t\leq m$. Then
$$V_{\mu+\Z e_{\bar i}+\cdots+\Z e_t}=\bigoplus_{\substack{k\in\Z_+\\i_1,\dots,i_k\in\{\bar i+1,\dots,t\}}}V_{\mu+e_{i_1}+\cdots+e_{i_k}-ke_{\bar i}},$$
where $V_{\mu+e_{i_1}+\cdots+e_{i_k}-ke_{\bar i}}=\C(\prod_{\substack{s=1}}^kE_{i_s,\bar i})\cdot v_\mu$.
\end{lemma}
\begin{proof}
Note that
$$V_{\mu+\Z e_{\bar i}+\cdots+\Z e_t}=\bigoplus_{\substack{k\in\Z_+\\j_1,\dots,j_k,j_1',\dots,j_k'\in\{\bar i,\dots,t\}\\j_1>j_1',\dots,j_k>j_k'}}\C E_{j_1,j_1'}\cdots E_{j_k,j_k'}\cdot v_\mu.$$
By the PBW theorem and Lemma \ref{Lem3.1}, we have $E_{j_1,j_1'}\cdots E_{j_k,j_k'}\cdot v_\mu=0$ when $\bar{i}+1\leq j_p'<j_p\leq t$ for any $1\leq p\leq k$. Therefore,
$$V_{\mu+\Z e_{\bar i}+\cdots+\Z e_t}=\bigoplus_{\substack{k\in\Z_+\\i_1,\dots,i_k\in\{\bar i+1,\dots,t\}}}\C E_{i_1,\bar i}\cdots E_{i_k,\bar i}\cdot v_\mu$$
and the lemma holds.
\end{proof}

\begin{lemma}\label{Lem3.3}
Suppose $\mu_{\bar i+1}=\dots=\mu_t$ for some $\bar i,t$ with $1\leq {\bar i}< t\leq m$ and $\mu_{\bar i}-\mu_{\bar i+1}=a$. Let $i_1,\dots,i_j, l\in\{\bar i+1,\dots,t\}$ satisfying $l\notin\{i_1,\ldots,i_j\}$ and $v=E_{i_1,\bar i}\cdots E_{i_j,\bar i}\cdot v_\mu\neq 0$.
\begin{itemize}
\item[(1)] If $a\in\N$, then $v$ generates a finite-dimensional simple module over the subalgebra generated by $E_{{\bar i},l}$ and $E_{l,\bar i}$, which is isomorphic to $\sl(2)$. Moreover, $E_{j_1,\bar i}\cdots E_{j_{a+1},\bar i}\cdot v_\mu=0$ for any $j_1,\ldots,j_{a+1}\in \{\bar i+1,\dots,t\}$.
\item[(2)] If $a\notin\Z_+$, then $v$ generates a infinite-dimensional simple module over the subalgebra generated by $E_{{\bar i},l}$ and $E_{l,\bar i}$, which is isomorphic to $\sl(2)$.
\end{itemize}
\end{lemma}
\begin{proof}
(1) Since $E_{\bar i,l}\cdot v\in V_{\mu+e_{i_1}+\dots+e_{i_j}-e_l-(j-1)e_{\bar i}}=0$ from Lemma \ref{Lem3.2} and $(E_{\bar i,\bar i}-E_{l,l})\cdot v=(a-j)v$, we obtain
$$E_{{\bar i},l}\cdot E_{l,\bar i}^{a-j+1}\cdot v=0$$
by the representation theory of $\sl(2)$. Let $w=E_{l,\bar i}^{a-j+1}\cdot v$. We claim that $E_{{\bar i}, r}\cdot w=0$ for any $r\in\{\bar i+1,\dots,t\}$. Of course, this is true when $r=l$.

If $r\in \{i_1,\ldots,i_j\}$, let's say $r=i_1$. Since $E_{{\bar i}, p}E_{{\bar i}, q}=E_{{\bar i}, q}E_{{\bar i}, p}$ for $p,q\in \{\bar i+1,\dots,t\}$, we can reorder $i_1,\dots,i_j$ such that $i_1=\cdots=i_n$ for $n<j$. Set $v'=E_{l,\bar i}^{a-j+1}E_{i_{n+1},\bar i}\cdots E_{i_j,\bar i}\cdot v_\mu$. We get $E_{{\bar i},i_1}\cdot E_{i_1,\bar i}^n\cdot v'=0$ since $E_{\bar i,i_1}\cdot v'=0$ and $(E_{\bar i,\bar i}-E_{i_1,i_1})\cdot v'=(n-1)v'$. Because
$$E_{i_1,\bar i}^n\cdot v'=E_{i_1,\bar i}^n\cdot E_{l,\bar i}^{a-j+1}E_{i_{n+1},\bar i}\cdots E_{i_j,\bar i}\cdot v_\mu=w,$$
we have $E_{\bar i,i_1}\cdot w=0$. So $E_{\bar i,r}\cdot w=0$ for $r\in \{i_1,\ldots,i_j\}$.

If $r\in\{\bar i+1,\dots,t\}\setminus \{i_1,\ldots,i_j,l\}$, we have
$$E_{\bar i,r}\cdot w\in V_{\mu+e_{i_1}+\cdots+e_{i_j}+(a-j+1)e_l-ae_{\bar i}-e_r}=0,$$
which implies that $E_{\bar i,r}\cdot w=0$ for $r\in \{\bar i+1,\dots,t\}\setminus \{i_1,\ldots,i_j,l\}$.

For any $1\leq i<r\leq m$ with $i\neq \bar i$ we have
$$E_{i,r}\cdot w\in V_{\mu+e_{i_1}+\cdots+e_{i_j}+(a-j+1)e_l+e_i-e_r-(a+1)e_{\bar i}}$$
and $E_{i,\bar i}\cdot w=0$. Then by the PBW theorem, $v_\mu\notin U(\gl(m)^+)\cdot w$. This contradicts the theory of simple highest weight modules. So $w=0$.

Therefore, the $\sl(2)$-module generated by $v$ as the highest weight vector has a basis $\{v, E_{l,\bar i}\cdot v,\ldots, E_{l,\bar i}^{a-j}\cdot v\}$ and is simple. Moreover, $w=0$ shows that $E_{j_1,\bar i}\cdots E_{j_{a+1},\bar i}\cdot v_\mu=0$ for any $j_1,\ldots,j_{a+1}\in \{\bar i+1,\dots,t\}$.

(2) Since $a\notin\Z_+$, we have $a-j\notin\Z_+$. Then the $\sl(2)$-module generated by $v$ as the highest weight vector has a basis
$\{v, E_{l,\bar i}\cdot v,\ldots, E_{l,\bar i}^{a-j+1}\cdot v,\ldots\}$ and is simple.
\end{proof}

By Lemma \ref{Lem3.3}, under the condition that $\mu_{\bar i+1}=\dots=\mu_t$ for some $\bar i$ with $1\leq {\bar i}< t\leq m$, we can make a more detailed characterization of the weight space $V_{\mu+\Z e_{\bar i}+\dots+\Z e_t}$ of $V$.

\begin{lemma}\label{Lem3.4}
Suppose $\mu_{\bar i+1}=\dots=\mu_t$ for some $\bar i$ with $1\leq {\bar i}< t\leq m$. Then
\begin{itemize}
\item[(1)] $$V_{\mu+\Z e_{\bar i}+\dots+\Z e_t}=\bigoplus_{\substack{k=0\\i_1,\dots,i_k\in\{\bar i+1,\dots,t\}}}^{\mu_{\bar i}-\mu_{\bar i+1}}V_{\mu+e_{i_1}+\dots+e_{i_k}-ke_{\bar i}}$$
if $\mu_{\bar i}-\mu_{\bar i+1}\in\N$, where $V_{\mu+e_{i_1}+\cdots+e_{i_k}-ke_{\bar i}}=\C (\prod_{\substack{s=1}}^kE_{i_s,\bar i})\cdot v_\mu\neq 0$ for $k\leqslant \mu_{\bar i}-\mu_{\bar i+1}$;
\item[(2)] $$V_{\mu+\Z e_{\bar i}+\cdots+\Z e_t}=\bigoplus_{\substack{k\in\Z_+\\i_1,\dots,i_k\in\{\bar i+1,\dots,t\}}}V_{\mu+e_{i_1}+\cdots+e_{i_k}-ke_{\bar i}}$$
if $\mu_{\bar i}-\mu_{\bar i+1}\notin\Z_+$, where $V_{\mu+e_{i_1}+\cdots+e_{i_k}-ke_{\bar i}}=\C (\prod_{\substack{s=1}}^kE_{i_s,\bar i})\cdot v_\mu\neq 0$.
\end{itemize}
\end{lemma}
\begin{proof}
By Lemma \ref{Lem3.2} and Lemma \ref{Lem3.3} (1), to prove (1), it suffices to prove that $(\prod_{\substack{s=1}}^kE_{i_s,\bar i})\cdot v_\mu\neq 0$ if and only if $0\leq k\leq \mu_{\bar i}-\mu_{\bar i+1}$. We will prove this by an induction on $k$ for $0\leq k\leqslant \mu_{\bar i}-\mu_{\bar i+1}$. If $k=0$, then $(\prod_{\substack{s=1}}^kE_{i_s,\bar i})\cdot v_\mu=v_\mu\neq 0$.

Reorder $i_1,\dots,i_k$ such that there is some $j\in\{1,\dots,k-1\}$ with $i_{j+1}=\cdots=i_k$ and $i_k\notin\{i_1,\dots,i_j\}$. By the induction hypothesis, we have $E_{i_1,\bar i}\cdots E_{i_j,\bar i}\cdot v_\mu\neq 0$. From Lemma \ref{Lem3.3} (1), we know that $E_{i_1,\bar i}\cdots E_{i_j,\bar i}\cdot v_\mu$ generates a simple $\sl(2)$-module, with a basis
$$E_{i_k,\bar i}^r\cdot E_{i_1,\bar i}\cdots E_{i_j,\bar i}\cdot v_\mu,\ 0\leq r\leq \mu_{\bar i}-\mu_{\bar i+1}-j.$$
This implies that $(\prod_{\substack{s=1}}^kE_{i_s,\bar i})\cdot v_\mu=E_{i_k,\bar i}^{k-j}\cdot E_{i_1,\bar i}\cdots E_{i_j,\bar i}\cdot v_\mu\neq 0$ since $0<k-j\leq \mu_{\bar i}-\mu_{\bar i+1}-j$.

The proof of (2) comes from Lemma \ref{Lem3.3} (2).
\end{proof}

For any $k\in\N$ and $l\in\Z$, set $A_l^{k}=l(l-1)\cdots(l-k+1)$ and $A_l^0=1$.
\begin{lemma}\label{Lem3.5}
Suppose $\mu_{\bar i+1}=\dots=\mu_t$ for some $\bar i$ with $1\leq {\bar i}< t\leq m$ and $\mu_{\bar i}\neq\mu_{\bar i+1}$.
Let $i_1,\dots,i_k\in\{\bar i+1,\dots,t\}$ and $l\in\N$ such that $\mu_{\bar i}-\mu_{\bar i+1}\in\C\setminus\{1,\dots,l-1\}$. Then there is a unique $v(i_1,\dots,i_k)\in V_{\mu+e_{i_1}+\dots+e_{i_k}-ke_{\bar i}}$ such that
$$(\prod_{\substack{s=1}}^kE_{\bar i,i_s})\cdot v(i_1,\dots,i_k)=A_l^kv_\mu.$$
\end{lemma}
\begin{proof}
By Lemma \ref{Lem3.2}, for any $i_1,\dots,i_k\in\{\bar i+1,\dots,t\}$, we have $V_{\mu+e_{i_1}+\dots+e_{i_k}-ke_{\bar i}}=\C E_{i_1,\bar i}\cdots E_{i_k,\bar i}\cdot v_\mu$. So $v(i_1,\dots,i_k)$ is unique if it exists. We will prove this lemma by an induction on $k$. Without loss of generality, suppose that $i_1\leqslant \dots\leqslant i_k$. For $k=1$, since
\begin{equation}\label{eq3.1}
E_{\bar i,i_1}\cdot (\mu_{\bar i}-\mu_{\bar i+1})^{-1}A_l^1E_{i_1,\bar i}\cdot v_\mu=A_l^1v_\mu,
\end{equation}
we get $E_{\bar i,i_1}\cdot v(i_1)=A_l^1v_\mu$ by taking $v(i_1)=(\mu_{\bar i}-\mu_{\bar i+1})^{-1}A_l^1E_{i_1,\bar i}\cdot v_\mu$.

Now suppose $k>1$ and for any $j<k$, we have
\begin{equation}\label{eq3.2}
(\prod_{\substack{s=1}}^jE_{\bar i,i_s})\cdot v(i_1,\dots,i_{j})=A_l^{j}v_\mu
\end{equation}
for a unique $v(i_1,\dots,i_{j})\in V_{\mu+e_{i_1}+\dots+e_{i_{j}}-je_{\bar i}}$. Let $j\in\{0,\dots,k-1\}$ with $i_j\neq i_{j+1}=\dots=i_k$. Then $v(i_1,\dots,i_{j})$ is a highest weight vector of a simple module over a copy of $\sl(2)$ that is generated by $E_{\bar i,i_k}$ and $E_{i_k,\bar i}$ from Lemma \ref{Lem3.3}. As we all know that this $\sl(2)$-module has a nonzero weight vector in $V_{\mu+e_{i_1}+\dots+e_{i_k}-ke_{\bar i}}$ if and only if $\mu_{\bar i}-\mu_{\bar i+1}-j\notin\Z_+$ or $k-j\leq \mu_{\bar i}-\mu_{\bar i+1}-j\in\Z_+$, in which case there exists a unique $v(i_1,\dots,i_k)\in V_{\mu+e_{i_1}+\dots+e_{\bar i+1}-ke_{\bar i}}$ such that
\begin{equation}\label{eq3.2-1}
(\prod_{\substack{s=j+1}}^kE_{\bar i,i_s})\cdot v(i_1,\dots,i_k)=A_{l-j}^{k-j}v(i_1,\dots,i_{j}).
\end{equation}
by the representation theory of $\sl(2)$. What remains is the case that $j\leq \mu_{\bar i}-\mu_{\bar i+1}\leq k-1$. In this case, $V_{\mu+e_{i_1}+\dots+e_{i_k}-ke_{\bar i}}=0$ by Lemma \ref{Lem3.4} (2). Also, $A_{l-j}^{k-j}=0$ since $k>\mu_{\bar i}-\mu_{\bar i+1}\geq l$. Let $v(i_1,\dots,i_k)=0$ and (\ref{eq3.2-1}) is still true.
Therefore, the combination \eqref{eq3.2} and \eqref{eq3.2-1} shows that there exists a unique $v(i_1,\dots,i_k)\in V_{\mu+e_{i_1}+\dots+e_{i_k}-ke_{\bar i}}$ such that
\begin{equation*}\begin{split}
(\prod_{\substack{s=1}}^kE_{\bar i,i_s})\cdot v(i_1,\dots,i_k)=&(\prod_{\substack{s=1}}^jE_{\bar i,i_s})\cdot(\prod_{\substack{p=j+1}}^kE_{\bar i,i_p})\cdot v(i_1,\dots,i_k)\\
=&A_{l-j}^{k-j}(\prod_{\substack{s=1}}^jE_{\bar i,i_s})\cdot v(i_1,\dots,i_{j})\\
=&A_{l}^kv_\mu.
\end{split}\end{equation*}
\end{proof}

Lemma \ref{Lem3.5} gives the relation between the weight vectors in $V_{\mu+e_{i_1}+\dots+e_{i_k}-ke_{\bar i}}$ and $v_\mu$ in the case $\mu_{\bar i}\neq\mu_{\bar i+1}$ and $\mu_{\bar i+1}=\dots=\mu_t$ for some $\bar i$ with $1\leq {\bar i}< t\leq m$. Using this lemma, we proceed to explore other relations between weight vectors in $L(\mu)$.

\begin{lemma}\label{Lem3.6}
Suppose $\mu_{\bar i+1}=\dots=\mu_t$ for some $\bar i$ with $1\leq {\bar i}< t\leq m$ and $\mu_{\bar i}\neq\mu_{\bar i+1}$.
Let $i,p,i_1,\dots,i_k\in\{\bar{i}+1,\dots, t\}$ with $i\neq p$, and $l\in\N$ such that $\mu_{\bar i}-\mu_{\bar i+1}\in\C\setminus\{1,\dots,l-1\}$. Then
\begin{itemize}
\item[(1)]$E_{i,\bar i}\cdot v_\mu=(\mu_{\bar i}-\mu_{\bar i+1})l^{-1}v(i)$;
\item[(2)]$E_{i,p}\cdot v(i_1,\dots,i_k)=0$ if $p\notin\{i_1,\dots,i_k\}$, and
$$E_{i,p}\cdot v(i_1,\dots,i_k)=(1+\sum_{s=1}^k\delta_{i_s,i})v(i_1,\dots,i_{j-1},i,i_{j+1},\dots,i_k)$$
if $p=i_j\in\{i_1,\dots,i_k\}$;
\item[(3)]$E_{i,\bar i}\cdot v(i_1,\dots,i_k)=\frac{1}{l-k}(\mu_{\bar i}-\mu_{\bar i+1}-k)(1+\sum_{s=1}^k\delta_{i_s,i})v(i_1,\dots,i_k,i)$
if $k<l$;
\item[(4)]$E_{\bar i,p}\cdot v(i_1,\dots,i_k)=0$ if $p\notin\{i_1,\dots,i_k\}$, and
$$E_{\bar i,p}\cdot v(i_1,\dots,i_k)=(l-k+1)v(i_1,\dots,i_{j-1},i_{j+1},\dots,i_k)$$
if $p=i_j\in\{i_1,\dots,i_k\}$.
\end{itemize}
\end{lemma}
\begin{proof}
Without loss of generality, assume $i_1\leqslant\dots\leqslant i_k$.

(1) From \eqref{eq3.1} we know $v(i)=(\mu_{\bar i}-\mu_{\bar i+1})^{-1}lE_{i,\bar i}\cdot v_\mu$. So $E_{i,\bar i}\cdot v_\mu=(\mu_{\bar i}-\mu_{\bar i+1})l^{-1}v(i)$.

(2) Note that $v(i_1,\dots,i_k)\in V_{\mu+e_{i_1}+\dots+e_{i_k}-ke_{\bar i}}$ and
$$E_{i,p}\cdot v(i_1,\dots,i_k)\in V_{\mu+e_{i_1}+\dots+e_{i_k}-ke_{\bar i}+e_i-e_p}.$$
If $p\notin\{i_1,\dots,i_k\}$, then we have $V_{\mu+e_{i_1}+\dots+e_{i_k}-ke_{\bar i}+e_i-e_p}=0$ by Lemma \ref{Lem3.2}.
So $E_{i,p}\cdot v(i_1,\dots,i_k)=0$. Now suppose that $p=i_j$ for some $j\in\{1,\dots,k\}$, then we have
\begin{equation*}\begin{split}
&E_{\bar i,i}(\prod_{\substack{s=1\\s\neq j}}^kE_{\bar i,i_s})\cdot E_{i,p}\cdot v(i_1,\dots,i_k)\\
=&\big(E_{i,p}E_{\bar i,i}(\prod_{\substack{s=1\\s\neq j}}^kE_{\bar i,i_s})+E_{\bar i,p}(\prod_{\substack{s=1\\s\neq j}}^kE_{\bar i,i_s})+E_{\bar i,i}[\prod_{\substack{s=1\\s\neq j}}^kE_{\bar i,i_s},E_{i,p}]\big)\cdot v(i_1,\dots,i_k)\\
=&0+(\prod_{\substack{s=1}}^k E_{\bar i,i_s})\cdot v(i_1,\dots,i_k)+(\sum_{\substack{s=1}}^k\delta_{i_s,i})(\prod_{\substack{s=1}}^k E_{\bar i,i_s})\cdot v(i_1,\dots,i_k)\\
=&(1+\sum_{s=1}^k\delta_{i_s,i})(\prod_{\substack{s=1}}^k E_{\bar i,i_s})\cdot v(i_1,\dots,i_k)\\
=&(1+\sum_{s=1}^k\delta_{i_s,i})A_l^k v_\mu.
\end{split}\end{equation*}
Lemma \ref{Lem3.3} shows that
$$E_{\bar i,i}(\prod_{\substack{s=1\\s\neq j}}^kE_{\bar i,i_s})\cdot v(i_1,\dots,i_{j-1},i,i_{j+1},\dots,i_k)=A_l^kv_\mu$$
and $ v(i_1,\dots,i_{j-1},i,i_{j+1},\dots,i_k)$ is unique. So we have
$$E_{i,p}\cdot v(i_1,\dots,i_k)=(1+\sum_{s=1}^k\delta_{i_s,i})v(i_1,\dots,i_{j-1},i,i_{j+1},\dots,i_k).$$

(3) Note that
\begin{equation*}\begin{split}
&E_{\bar i,i}(\prod_{\substack{s=1}}^kE_{\bar i,i_s})\cdot E_{i,\bar i}\cdot v(i_1,\dots,i_k)\\
=&E_{\bar i,i}\cdot\big([\prod_{\substack{s=1}}^kE_{\bar i,i_s}, E_{i,\bar i}]+E_{i,\bar i}(\prod_{\substack{s=1}}^kE_{\bar i,i_s})\big)\cdot v(i_1,\dots,i_k)\\
=&E_{\bar i,i}\big(\sum_{s=1}^kE_{\bar i,i_1}\cdots E_{\bar i,i_{s-1}}(\delta_{i_s,i}E_{\bar i,\bar i}-E_{i,i_s})E_{\bar i,i_{s+1}}\cdots E_{\bar i,i_k}\big)\cdot v(i_1,\dots,i_k)+E_{\bar i,i}E_{i,\bar i}\cdot A_l^kv_\mu\\
=&(\sum_{s=1}^k\delta_{i_s,i}(\mu_{\bar i}-s))(\prod_{\substack{s=1}}^kE_{\bar i,i_s})\cdot v(i_1,\dots,i_k)+(\sum_{s=1}^k\delta_{i_{s+1},i}+\dots+\delta_{i_k,i})(\prod_{\substack{s=1}}^kE_{\bar i,i_s})\cdot v(i_1,\dots,i_k)
\end{split}\end{equation*}
\begin{equation*}\begin{split}
&-E_{\bar i,i}(\sum_{s=1}^kE_{\bar i,i_1}\cdots E_{\bar i,i_{s-1}}E_{\bar i,i_{s+1}}\cdots E_{\bar i,i_k}E_{i,i_s})\cdot v(i_1,\dots,i_k)+(\mu_{\bar i}-\mu_{\bar i+1})A_l^kv_\mu\\
=&(\sum_{s=1}^k\delta_{i_s,i}(\mu_{\bar i}-s))A_l^kv_\mu+\sum_{s=1}^k(\delta_{i_{s+1},i}+\dots+\delta_{i_k,i})A_l^kv_\mu\\
&-E_{\bar i,i}\big(\sum_{s=1}^kE_{\bar i,i_1}\cdots E_{\bar i,i_{s-1}}E_{\bar i,i_{s+1}}\cdots E_{\bar i,i_k}
(1-\delta_{i_s,i})(1+\sum_{r=1}^k\delta_{i_r,i})\\
&\cdot v(i_1,\dots,i_{s-1},i,i_{s+1},\dots,i_k)\big)\\
&-E_{\bar i,i}(\sum_{s=1}^kE_{\bar i,i_1}\cdots E_{\bar i,i_{s-1}}E_{\bar i,i_{s+1}}\cdots E_{\bar i,i_k}\delta_{i_s,i})(\mu_{\bar i+1}+\sum_{r=1}^k\delta_{i_r,i})\cdot v(i_1,\dots,i_k)+(\mu_{\bar i}-\mu_{\bar i+1})A_l^kv_\mu\\
=&(\mu_{\bar i}-1)(\sum_{s=1}^k\delta_{i_s,i})A_l^kv_\mu-(\sum_{s=1}^k(1-\delta_{i_s,i}))(1+\sum_{r=1}^k\delta_{i_r,i})A_l^kv_\mu\\
&-(\sum_{s=1}^k\delta_{i_s,i})(\mu_{\bar i+1}+\sum_{r=1}^k\delta_{i_r,i})A_l^kv_\mu+(\mu_{\bar i}-\mu_{\bar i+1})A_l^kv_\mu\\
=&(\mu_{\bar i}-\mu_{\bar i+1}-k)(1+\sum_{s=1}^k\delta_{i_s,i})A_l^kv_\mu.
\end{split}\end{equation*}
Then, we have
\begin{equation}\label{eq3.3}
(l-k)E_{\bar i,i}(\prod_{\substack{s=1}}^kE_{\bar i,i_s})\cdot E_{i,\bar i}\cdot v(i_1,\dots,i_k)=(\mu_{\bar i}-\mu_{\bar i+1}-k)(1+\sum_{s=1}^k\delta_{i_s,i})A_l^{k+1}v_\mu.
\end{equation}
By Lemma \ref{Lem3.5}, we have
\begin{equation}\label{eq3.4}
E_{i,\bar i}\cdot (\prod_{\substack{s=1}}^kE_{\bar i,i_s})\cdot E_{i,\bar i}\cdot v(i_1,\dots,i_k,i)=A_l^{k+1}v_\mu.
\end{equation}
Since $l\leq \mu_{\bar i}-\mu_{\bar i+1}$, we have $\mu_{\bar i}-\mu_{\bar i+1}-k\neq 0$ when $k<l$. Compare \eqref{eq3.3} and \eqref{eq3.4}, and we get
$$(l-k)E_{i,\bar i}\cdot v(i_1,\dots,i_k)=(\mu_{\bar i}-\mu_{\bar i+1}-k)(1+\sum_{s=1}^k\delta_{i_s,i})v(i_1,\dots,i_k,i)$$
if $k<l$.

(4) If $p\notin\{i_1,\dots,i_k\}$, then $E_{\bar i,p}\cdot v(i_1,\dots,i_k)=0$ by Lemma \ref{Lem3.3} (1). Now suppose $p=i_j$ for some $j\in\{1,\dots,k\}$, then
$$(\prod_{\substack{s=1\\s\neq j}}^kE_{\bar i,i_s})\cdot E_{\bar i,p}\cdot v(i_1,\dots,i_k)=A_l^kv_\mu=(l-k+1)A_l^{k-1}v_\mu$$
by Lemma \ref{Lem3.5}. Then (4) is also true.
\end{proof}

\section{Irreducibility of $M_{\p}(V)$}
It is mentioned in Subsection 2.2 that the subalgebra $\mathfrak l$ of $\sl(m+1)$ is isomorphic to $\gl(m)$. In this section, let $V$ be the simple highest weight $\mathfrak l$-module $L(\mu)$ with highest weight vector $v_\mu$ for some $\mu\in\C^m$, where $(E_{i,i}-E_{m+1,m+1})\cdot v_\mu=\mu_iv_\mu$. The results in Section 3 concerning the simple highest weight module of $\gl(m)$ are equally applicable to the $\mathfrak l$-module $L(\mu)$.

Let $\mathfrak{u}^-=\sum\limits_{s=1}^m\C E_{m+1,s}$ and
$$Y_\alpha=E_{m+1,1}^{\alpha_1}E_{m+1,2}^{\alpha_2}\cdots E_{m+1,m}^{\alpha_m}\in U(\mathfrak{u}^-)$$
for any $\alpha\in \Z_+^m$. Set $Y_\alpha=1\in U(\mathfrak{u}^-)$ for any $\alpha\in\Z^m\setminus \Z_+^m$.
Since $\s=\p\oplus \mathfrak{u}^-$, we know that $M_{\p}(V)\simeq U(\mathfrak{u}^-)\otimes_\C V$ as $U(\mathfrak{u}^-)$-modules by the PBW theorem. Therefore, for $\alpha\in\Z_+^m$ and $v\in V$, we have $Y_\alpha v\in M_{\p}(V)$ and the following computations.

\begin{lemma}\label{Lem4.1}
Let $\alpha\in\Z_+^m$ and $v\in V$. Then
\begin{itemize}
\item[(1)] $E_{i,j}\cdot(Y_\alpha v)=-\alpha_iY_{\alpha-e_i+e_j} v+Y_\alpha E_{i,j}\cdot v$;
\item[(2)] $(E_{i,i}-E_{j,j})\cdot(Y_\alpha v)=(\alpha_j-\alpha_i)Y_\alpha v+Y_\alpha(E_{i,i}-E_{j,j})\cdot v;$
\item[(3)] $(E_{i,i}-E_{m+1,m+1})\cdot (Y_\alpha v)=-(\alpha_i+|\alpha|)Y_\alpha v+Y_\alpha(E_{i,i}-E_{m+1,m+1})\cdot v;$
\item[(4)] $E_{m+1,i}\cdot (Y_\alpha v)=Y_{\alpha+e_i} v;$
\item[(5)] $E_{i,m+1}\cdot (Y_\alpha v)=\alpha_i(1-|\alpha|)Y_{\alpha-e_i} v+\sum\limits_{\substack{j=1\\j\neq i}}^m\alpha_jY_{\alpha-e_j} E_{i,j}\cdot v+\alpha_iY_{\alpha-e_i}(E_{i,i}-E_{m+1,m+1})\cdot v$
\end{itemize}
for $1\leq i\neq j\leq m$.
\end{lemma}
\begin{proof}
(1)-(4) can be obtained by direct verification.

(5)
\begin{eqnarray}\label{eq4.1}
&&E_{i,m+1}\cdot(Y_\alpha v)\\
&=&\sum_{j=1}^m\sum_{s=1}^{\alpha_j}Y_{\alpha_1e_1+\dots+\alpha_{j-1}e_{j-1}+(s-1)e_j}(E_{i,j}-\delta_{i,j}E_{m+1,m+1})\cdot Y_{(\alpha_j-s)e_j+\alpha_{j+1}e_{j+1}+\dots+\alpha_me_m}v\nonumber
\end{eqnarray}
By (1) and (3), we have
\begin{eqnarray*}
&&E_{i,j}\cdot Y_{(\alpha_j-s)e_j+\alpha_{j+1}e_{j+1}+\dots+\alpha_me_m}v=Y_{(\alpha_j-s)e_j+\alpha_{j+1}e_{j+1}+\dots+\alpha_me_m}E_{i,j}\cdot v\\
&&-(\alpha_{j+1}e_{j+1}+\dots+\alpha_me_m)_iY_{-e_i+e_j+(\alpha_j-s)e_j+\alpha_{j+1}e_{j+1}+\dots+\alpha_me_m}v,
\end{eqnarray*}
\begin{eqnarray*}
&&(E_{i,i}-E_{m+1,m+1})\cdot Y_{(\alpha_j-s)e_j+\alpha_{j+1}e_{j+1}+\dots+\alpha_me_m}v\\
&=&Y_{(\alpha_j-s)e_j+\alpha_{j+1}e_{j+1}+\dots+\alpha_me_m}(E_{i,i}-E_{m+1,m+1})\cdot v\\
&&-(\alpha_i-s+\alpha_i-s+\alpha_{i+1}+\dots+\alpha_m)Y_{(\alpha_j-s)e_j+\alpha_{j+1}e_{j+1}+\dots+\alpha_me_m}v.
\end{eqnarray*}
By substituting the above two equations into (\ref{eq4.1}), we obtain
\begin{eqnarray*}
&&E_{i,m+1}\cdot(Y_\alpha v)\\
&=&\sum_{\substack{j=1\\j\neq i}}^m\sum_{s=1}^{\alpha_j}\big(-(\alpha_{j+1}e_{j+1}+\dots+\alpha_me_m)_iY_{\alpha-e_i}v+Y_{\alpha-e_j}E_{i,j}\cdot v\big)\\
&&+\sum_{s=1}^{\alpha_i}\big(-(2\alpha_i-2s+\alpha_{i+1}+\dots+\alpha_m)Y_{\alpha-e_i}v+Y_{\alpha-e_i}(E_{i,i}-E_{m+1,m+1})\cdot v\big)\\
&=&-\alpha_i(\alpha_1+\dots+\alpha_{i-1})Y_{\alpha-e_i}v+\sum_{\substack{j=1\\j\neq i}}^m\alpha_jY_{\alpha-e_j}E_{i,j}\cdot v\\
&&+\big(-\alpha_i(\alpha_i+\alpha_{i+1}+\dots+\alpha_m)+\alpha_i\big)Y_{\alpha-e_i}v+\alpha_iY_{\alpha-e_i}(E_{i,i}-E_{m+1,m+1})\cdot v\\
&=&\alpha_i(1-|\alpha|)Y_{\alpha-e_i} v+\sum\limits_{\substack{j=1\\j\neq i}}^m\alpha_jY_{\alpha-e_j}E_{ij}\cdot v+\alpha_iY_{\alpha-e_i}(E_{i,i}-E_{m+1,m+1})\cdot v.
\end{eqnarray*}
\end{proof}

We define the strict total order ``$\prec$" in $\Z_+^m$ as follows. For any $\alpha, \beta\in \Z_+^m$, we say that $\alpha\prec\beta$ if $|\alpha|<|\beta|$ or
$|\alpha|=|\beta|, \alpha_1=\beta_1,\dots, \alpha_{i-1}=\beta_{i-1}, \alpha_i<\beta_i$ for some $i\in\{1,\dots,m\}$.
For any $u\in M_{\p}(V)$, $u$ is uniquely of the form $\sum_{\alpha\in\Z_+^m}Y_\alpha\cdot u(\alpha)$, where $u(\alpha)\in V$. Set
$$S_u=\{\alpha\in\Z_+^m\ |\ u(\alpha)\neq 0\}$$
and $u(\alpha)=0$ for all $\alpha\in\Z^m\setminus\Z_+^m$. For $u\neq 0$, let $\alpha_u$ be the maximal element in $S_u$ with respect to ``$\succ$". We call $Y_{\alpha_u}\cdot u(\alpha_u)$ the \emph{leading term} of $u$.

Now suppose that $M_{\p}(V)$ is not a simple $\s$-module. Let $w=\sum_{\alpha\in\Z_+^m}Y_\alpha\cdot w(\alpha)\in M_{\p}(V)\setminus\{\C v_\mu\}$ be a highest weight vector with minimal $\alpha_w$. Next, we will discuss some properties of $w$.
\begin{lemma}\label{Lem4.2}
The vector $w$ is unique up to a nonzero scalar and $w(\alpha_w)\in\C v_\mu$.
\end{lemma}
\begin{proof}
If $w(\alpha_w)\notin\C v_\mu$, then there exist $i,j\in\{1,\dots,m\}$ with $i<j$ such that $E_{i,j}\cdot w(\alpha_w)\neq 0$. By Lemma \ref{Lem4.1} (1), $(E_{i,j}\cdot w)(\alpha_w)=E_{i,j}\cdot w(\alpha_w)\neq 0$. This contradicts the fact that $w$ is a highest weight vector. So $w(\alpha_w)\in\C v_\mu$. Clearly, such $w$ is unique up to a nonzero scalar.
\end{proof}
Hence, there is a unique highest weight vector $w\in M_{\p}(V)\setminus\{\C v_\mu\}$ such that $\alpha_w$ is minimal and the leading term of $w$ is $Y_{\alpha_w}v_\mu$. Fix this $w$ and denote $\alpha_w$ by $\bar\alpha$.

\begin{lemma}\label{Lem4.3}
\begin{itemize}
\item[(1)] For $1\leq i\leq m$, we have $(E_{i,i}-E_{m+1,m+1})\cdot w=(-\bar\alpha_i-|\bar\alpha|+\mu_i)w$.
\item[(2)] For any $\beta\in\Z^m$, we have $w(\bar\alpha+\beta)\in V_{\mu+\beta}$.
\end{itemize}
\end{lemma}
\begin{proof}
By Lemma \ref{Lem4.1} (3), we have
$$(E_{i,i}-E_{m+1,m+1})\cdot (Y_\alpha w(\alpha))=-(\alpha_i+|\alpha|)Y_\alpha w(\alpha)+Y_\alpha(E_{i,i}-E_{m+1,m+1})\cdot w(\alpha)$$
for any $\alpha\in S_w$. In particular,
$$(E_{i,i}-E_{m+1,m+1})\cdot (Y_{\bar{\alpha}}v_\mu)=(-\bar\alpha_i-|\bar\alpha|+\mu_i)Y_{\bar{\alpha}} v_\mu.$$
Since $w$ is a weight vector, we have $(E_{i,i}-E_{m+1,m+1})\cdot w=(-\bar\alpha_i-|\bar\alpha|+\mu_i)w$.
Consequently,
$$-(\alpha_i+|\alpha|)w(\alpha)+(E_{i,i}-E_{m+1,m+1})\cdot w(\alpha)=(-\bar\alpha_i-|\bar\alpha|+\mu_i)w(\alpha)$$
for any $\alpha\in S_w$. Then for any $\beta\in\Z^m$,
$$(E_{i,i}-E_{m+1,m+1})\cdot w(\bar\alpha+\beta)=(\beta_i+|\beta|+\mu_i)w(\bar\alpha+\beta),$$
which implies that $w(\bar\alpha+\beta)\in V_{\mu+\beta+|\beta|(e_1+\dots+e_m)}$. Note that $V=\bigoplus_{\substack{\beta\in\Z^m\\|\beta|=0}}V_{\mu+\beta}$. If $|\beta|\neq 0$, then $\big|\beta+|\beta|(e_1+\dots+e_m)\big|=(m+1)|\beta|\neq 0$, in which case $w(\bar\alpha+\beta)\in V_{\mu+\beta+|\beta|(e_1+\dots+e_m)}=0$. So $w(\bar\alpha+\beta)\in V_{\mu+\beta}$ for any $\beta\in\Z^m$.
\end{proof}

\begin{lemma}\label{Lem4.4}
For any $\beta\in\Z^m$ and $1\leq i< j\leq m$, we have
$$E_{i,j}\cdot w(\bar\alpha+\beta)=(\bar\alpha_i+\beta_i+1) w(\bar\alpha+\beta+e_i-e_j).$$
In particular, $E_{i,j}\cdot w(\bar\alpha-e_i+e_j)=\bar\alpha_iv_\mu$.
\end{lemma}
\begin{proof}
Let $w_1=E_{i,j}\cdot w$. Then $w_1=0$ by the definition of $w$, i.e., $w_1(\bar\alpha+\beta)=0$ for any $\beta\in\Z^m$. Meanwhile, if $\bar\alpha+\beta\in S_{E_{i,j}\cdot (Y_\alpha\cdot w(\alpha))}$, then $\alpha=\bar\alpha+\beta$ or $\alpha=\bar\alpha+\beta+e_i-e_j$ by Lemma \ref{Lem4.1} (1). Since $w_1(\bar\alpha+\beta)=0$ we have
$$E_{i,j}\cdot w(\bar\alpha+\beta)-(\bar\alpha_i+\beta_i+1)w(\bar\alpha+\beta+e_i-e_j)=0.$$
In particular, if $\beta=-e_i+e_j$, then $E_{i,j}\cdot w(\bar\alpha-e_i+e_j)=\bar\alpha_iw(\bar\alpha)=\bar\alpha_iv_\mu$.
\end{proof}
\begin{corollary}\label{Cor4.5}
If $\bar\alpha_j\neq 0$ for some $j\in\{1,\dots,m-1\}$, then $\mu_j-\mu_{j+1}\in\C\setminus\{0,\dots,\bar\alpha_j-1\}$.
\end{corollary}
\begin{proof}
By Lemma \ref{Lem4.4}, we can get $$E_{j,j+1}^k\cdot w(\bar\alpha+k(e_{j+1}-e_j))=A_{\bar\alpha_j}^kv_\mu,\ \forall k\in\Z_+.$$ So $ w(\bar\alpha+k(e_{j+1}-e_j))\neq 0$ for $k=0,\dots,\bar\alpha_j$. By Lemma \ref{Lem4.3} (2), $w(\bar\alpha+k(e_{j+1}-e_j))\in V_{\mu+k(e_{j+1}-e_j)}$. Then $\oplus_{k\in\Z_+} V_{\mu+k(e_{j+1}-e_j)}\geqslant\bar\alpha_i+1$.

Note that the space $\oplus_{k\in\Z_+} V_{\mu+k(e_{j+1}-e_j)}$ is a simple highest weight module over the subalgebra of $\gl(m)$ that is generated by $E_{j,j+1}$ and $E_{j+1,j}$, which is isomorphic to $\sl(2)$. Suppose $\mu_j-\mu_{j+1}\in\Z_+$. Then $\bar\alpha_i+1\leqslant\dim\bigoplus_{k\in\Z_+} V_{\mu+k(e_{j+1}-e_j)}=\mu_j-\mu_{j+1}+1$. So $\bar\alpha_j\leqslant\mu_j-\mu_{j+1}$.
\end{proof}

\begin{lemma}\label{Lem4.6}
Suppose $1\leq i\leq m$. Then we have
\begin{itemize}
\item[(1)]for any $\beta\in\Z^m$,
\begin{equation*}\begin{split}
&\big(\sum_{s=1}^i(\bar\alpha_s+1)-\sum_{s=i}^m \beta_s-\mu_i-m-1\big)(\bar\alpha_i+\beta_i)w(\bar\alpha+\beta)\\
=&\sum_{s=1}^{i-1}(\bar\alpha_{s}+\beta_{s}+1)E_{i,s}\cdot w(\bar\alpha+\beta-e_i+e_{s});
\end{split}\end{equation*}
\item[(2)] $\big(\sum_{s=1}^i(\bar\alpha_s+1)-\mu_i-m-1\big)\bar\alpha_i=0.$
\end{itemize}
\end{lemma}
\begin{proof}
(1) If $|\beta|\neq 0$, then $w(\bar\alpha+\beta)=w(\bar\alpha+\beta-e_i+e_s)=0$ for $1\leq s\leq i-1$ by Lemma \ref{Lem4.3} (2). Then (1) holds. Now assume $|\beta|=0$.

By Lemma \ref{Lem4.1} (5), if $\bar\alpha+\beta-e_i\in S_{E_{i,m+1}\cdot(Y_\alpha\cdot w(\alpha))}$, then $\alpha=\bar\alpha+\beta-e_i+e_s$ for some $s\in\{1,\dots,m\}$. From Lemma \ref{Lem4.1} (5) and \ref{Lem4.3} (2), we have
\begin{equation*}\begin{split}
&{\big(E_{i,m+1}\cdot (Y_{\bar\alpha+\beta} w(\bar\alpha+\beta))\big)}(\bar\alpha+\beta-e_i)\\
=&(\bar\alpha_i+\beta_i)(1-|\bar\alpha|)w(\bar\alpha+\beta)+(\bar\alpha_i+\beta_i)E_{i,i}\cdot w(\bar\alpha+\beta)\\
=&(\bar\alpha_i+\beta_i)(1-|\bar\alpha|+\mu_i+\beta_i)w(\bar\alpha+\beta).
\end{split}\end{equation*}
and
\begin{equation*}\begin{split}
&{\big(E_{i,m+1}\cdot (Y_{\bar\alpha+\beta-e_i+e_s} w(\bar\alpha+\beta-e_i+e_s))\big)}(\bar\alpha+\beta-e_i)\\
=&(\bar\alpha_s+\beta_s+1)E_{i,s}\cdot w(\bar\alpha+\beta-e_i+e_s)
\end{split}\end{equation*}
for $1\leq i\neq s\leq m$. By the definition of $w$, we know that $E_{i,m+1}\cdot w=0$ for $1\leq i\leq m$, which implies that ${(E_{i,m+1}\cdot w)}(\bar\alpha+\beta-e_i)=0$ for any $\beta\in\Z^m$. Therefore, we get
$$(\bar\alpha_i+\beta_i)(|\bar\alpha|-\mu_i-\beta_i-1)w(\bar\alpha+\beta)$$
$$=\sum_{\substack{s=1\\s\neq i}}^m(\bar\alpha_s+\beta_s+1)E_{i,s}\cdot w(\bar\alpha+\beta-e_i+e_s).$$
By Lemma \ref{Lem4.4}, $(\bar\alpha_{s}+\beta_{s}+1)E_{i,s}\cdot w(\bar\alpha+\beta-e_i+e_{s})=(\bar\alpha_{s}+\beta_{s}+1)(\bar\alpha_i+\beta_i) w(\bar\alpha+\beta)$ if $s>i$. Then
\begin{equation*}\begin{split}
&\big(\sum_{s=1}^i(\bar\alpha_s+1)-\sum_{s=i}^m \beta_s-\mu_i-m-1\big)(\bar\alpha_i+\beta_i)w(\bar\alpha+\beta)\\
=&\sum_{s=1}^{i-1}(\bar\alpha_{s}+\beta_{s}+1)E_{i,s}\cdot w(\bar\alpha+\beta-e_i+e_{s}).
\end{split}\end{equation*}
(2) Suppose $\beta=\bf 0$, then $w(\bar\alpha-e_i+e_{s})=0$ if $s<i$ by the definition of $w$. From (1), we have
$$\big(\sum_{s=1}^i(\bar\alpha_s+1)-\mu_i-m-1\big)\bar\alpha_i=0.$$
\end{proof}

\begin{corollary}\label{Cor4.7}
Suppose that $V$ is finite-dimensional. Then there is some $i\in\{1,\dots,m\}$ such that $\mu_i+m+1-i\in\N$ and $\bar\alpha=(\mu_i+m+1-i)e_i$.
\end{corollary}
\begin{proof}
Suppose that there is $i,j\in\{1,\dots,m\}$ with $i<j$ such that $\bar\alpha_i,\bar\alpha_j\neq 0$. By Lemma \ref{Lem4.6} (2),
$$\sum_{s=1}^i(\bar\alpha_s+1)-(\mu_i+m+1)=\sum_{s=1}^j(\bar\alpha_s+1)-(\mu_j+m+1)=0.$$
Then $\mu_j-\mu_i=\sum_{s=i+1}^j(\bar\alpha_s+1)>0$. Since $V=L(\mu)$ is finite-dimensional, we have $\mu_1\geqslant\dots\geqslant\mu_m$, which is a contradiction.
So $\bar\alpha=\bar\alpha_ie_i$ for some $i$ and $\bar\alpha_i=\mu_i+m+1-i\in\N$.
\end{proof}

\begin{corollary}
Suppose $\mu_i+m+1\notin\N$ for all $i\in\{1,\dots,m\}$. Then $M_{\p}(V)$ is a simple $\s$-module.
\end{corollary}
\begin{proof}
Suppose that $M_{\p}(V)$ is not simple. Let $i$ be the smallest in $\{1,\dots,m\}$ such that $\bar\alpha_i\neq 0$. By Lemma \ref{Lem4.6} (2), $\mu_i+m+1=\bar{\alpha}_i+i\in\N$, which is a contradiction.
\end{proof}

We also need the following lemma on highest weight vectors.
\begin{lemma}\label{Lem4.9}
Suppose that $\mu_m+1\in\N$. Then $Y_{(\mu_m+1)e_m}v_\mu$ is a highest weight vector of $M_{\p}(V)$.
\end{lemma}
\begin{proof}
By Lemma \ref{Lem4.1} (1), for any $i,j\in\{1,\dots,m\}$ with $i<j$, we have $E_{i,j}\cdot Y_{(\mu_m+1)e_m}v_\mu=0$. Also, by Lemma \ref{Lem4.1} (5), we obtain
\begin{equation*}\begin{split}
&E_{m,m+1}\cdot Y_{(\mu_m+1)e_m}v_\mu\\
=&(\mu_m+1)(1-\mu_m-1)Y_{\mu_me_m}v_\mu+(\mu_m+1)Y_{\mu_me_m}(E_{m,m}-E_{m+1,m+1})\cdot v_\mu=0.
\end{split}\end{equation*}
So $Y_{(\mu_m+1)e_m}v_\mu$ is a highest weight vector of $M_{\p}(V)$.
\end{proof}

In the rest of this section, we will focus on the case that
$$\mu_1=\dots=\mu_{\bar i},\ \mu_{\bar i+1}=\dots=\mu_m$$
for some $\bar i\in\{1,\dots,m\}$. Set $\bar i=m$ if $\mu_1=\dots=\mu_m$.
We will give the necessary and sufficient conditions for $M_{\p}(V)$ not to be a simple $\s$-module.

For the case that $\bar i\neq m$, let $l=\mu_{\bar i}+m+1-\bar i$. If $l\in\N$ and $\mu_{\bar i}-\mu_{\bar i+1}\in\C\setminus\{0,\dots,l-1\}$, let $v(i_1,\dots,i_k)\in V_{\mu+e_{i_1}+\dots+e_{i_k}-ke_{\bar i}}$ be as in Lemma \ref{Lem3.5} for any $i_1,\dots,i_k\in\{\bar i+1,\dots,m\}$. Define the following vector $u$ in $M_\p(V)$:
$$u=\sum_{k=0}^l\sum_{\substack{i_1,\dots,i_k\in\{\bar i+1,\dots,m\}\\i_1\leqslant\dots\leqslant i_k}}Y_{le_{\bar i}+e_{i_1}+\dots+e_{i_k}-ke_{\bar i}}v(i_1,\dots,i_k).$$
\begin{lemma}\label{Lem4.10}
Suppose $\bar i\neq m,l\in\N$ and $\mu_{\bar i}-\mu_{\bar i+1}\in\C\setminus\{0,\dots,l-1\}$. Then $u$ is a highest weight vector of $M_\p(V)$.
\end{lemma}
\begin{proof}
To show that $u$ is a highest weight vector, it suffices to verify that $E_{i,j}\cdot u=0$ for any $i,j\in\{1,\dots,m\}$ with $i<j$ and $E_{m,m+1}\cdot u=0$.

\textbf{Claim 1.} $E_{i,j}\cdot u=0$.

Let $i,j\in\{1,\dots,m\}$ with $i<j$. For any $i'_1,\dots,i'_{k'}\in\{\bar i+1,\dots,m\}$,
\begin{equation*}\begin{split}
&E_{i,j}\cdot Y_{le_{\bar i}+e_{i'_1}+\dots+e_{i'_{k'}}-k'e_{\bar i}}v(i'_1,\dots,i'_{k'})\\
=&-(le_{\bar i}+e_{i'_1}+\dots+e_{i'_{k'}}-k'e_{\bar i})_iY_{le_{\bar i}+e_{i'_1}+\dots+e_{i'_{k'}}-k'e_{\bar i}-e_i+e_j}v(i'_1,\dots,i'_{k'})\\
&+Y_{le_{\bar i}+e_{i'_1}+\dots+e_{i'_{k'}}-k'e_{\bar i}}E_{i,j}\cdot v(i'_1,\dots,i'_{k'})
\end{split}\end{equation*}
by Lemma \ref{Lem4.1} (1).

If $i<\bar i$, then $E_{i,j}\cdot v(i'_1,\dots,i'_{k'})=0$ by the structure of the highest weight module $V$, and $(le_{\bar i}+e_{i'_1}+\dots+e_{i'_{k'}}-k'e_{\bar i})_i=0$. It follows that $E_{i,j}\cdot Y_{le_{\bar i}+e_{i'_1}+\dots+e_{i'_{k'}}-k'e_{\bar i}}v(i'_1,\dots,i'_{k'})=0$. So $E_{i,j}\cdot u=0$.

Now suppose that $i\geqslant\bar i$. Note that any element in $S_{E_{i,j}\cdot Y_{le_{\bar i}+e_{i'_1}+\dots+e_{i'_{k'}}-k'e_{\bar i}}v(i'_1,\dots,i'_{k'})}$ is of the form $le_{\bar i}+e_{i_1}+\dots+e_{i_k}-ke_{\bar i}$ for some $i_1,\dots,i_k\in\{\bar i+1,\dots,m\}$. Moreover, if $le_{\bar i}+e_{i_1}+\dots+e_{i_k}-ke_{\bar i}\in S_{E_{i,j}\cdot Y_{le_{\bar i}+e_{i'_1}+\dots+e_{i'_{k'}}-k'e_{\bar i}}v(i'_1,\dots,i'_{k'})}$, then one of the following is true:
\begin{itemize}
\item[(1)] $e_{i'_1}+\dots+e_{i'_{k'}}=e_{i_1}+\dots+e_{i_k}$;
\item[(2)] $e_{i'_1}+\dots+e_{i'_{k'}}=e_{i_1}+\dots+e_{i_k}-e_j+e_i,\ i>\bar i,j\in\{i_1,\dots,i_k\},k'=k>0$;
\item[(3)] $e_{i'_1}+\dots+e_{i'_{k'}}=e_{i_1}+\dots+e_{i_k}-e_j,\ i=\bar i,j\in\{i_1,\dots,i_k\},k'+1=k>0$.
\end{itemize}
If $j\notin\{i_1,\dots,i_k\}$, then $E_{i,j}\cdot v(i_1,\dots,i_k)=0$ by Lemma \ref{Lem3.6} (2), (4). If $k=0$, then $E_{i,j}\cdot v(i_1,\dots,i_k)=E_{i,j}\cdot v_\mu=0$. So
$$\big(E_{i,j}\cdot u\big)(le_{\bar i}+e_{i_1}+\dots+e_{i_k}-ke_{\bar i})=0$$
if $j\notin\{i_1,\dots,i_k\}$ or $k=0$.

Now consider the case that $k>0$ and $j=i_r$ for some $r\in\{1,\dots,k\}$. Then
\begin{equation*}\begin{split}
&\big(E_{i,j}\cdot Y_{le_{\bar i}+e_{i_1}+\dots+e_{i_k}-ke_{\bar i}}v(i_1,\dots,i_k)\big)(le_{\bar i}+e_{i_1}+\dots+e_{i_k}-ke_{\bar i})\\
=&E_{i,j}\cdot v(i_1,\dots,i_k)\\
=&\left\{\begin{aligned}&(1+\sum_{n=1}^k\delta_{i_n,i})v(i_1,\dots,i_{r-1},i,i_{r+1},\dots,i_k), &\text{ if }i>\bar i;\\&(l-k+1)v(i_1,\dots,i_{r-1},i_{r+1},\dots,i_k),&\text{ if }i=\bar i,\end{aligned}\right.
\end{split}\end{equation*}
from Lemma \ref{Lem3.6}.
If $i>\bar i$, then
\begin{equation*}\begin{split}
&\big(E_{i,j}\cdot Y_{le_{\bar i}+e_{i_1}+\dots+e_{i_k}-e_j+e_i-ke_{\bar i}}v(i_1,\dots,i_{r-i},i,i_{r+1},\dots,i_k)\big)(le_{\bar i}+e_{i_1}+\dots+e_{i_k}-ke_{\bar i})\\
=&-(1+\sum_{n=1}^k\delta_{i_n,i})v(i_1,\dots,i_{r-1},i,i_{r+1},\dots,i_k).
\end{split}\end{equation*}
If $i=\bar i$, then
\begin{equation*}\begin{split}
&\big(E_{i,j}\cdot Y_{le_{\bar i}+e_{i_1}+\dots+e_{i_k}-e_j-(k-1)e_{\bar i}}v(i_1,\dots,i_{r-i},i_{r+1},\dots,i_k)\big)(le_{\bar i}+e_{i_1}+\dots+e_{i_k}-ke_{\bar i})\\
=&-(l-k+1)v(i_1,\dots,i_{r-1},i,i_{r+1},\dots,i_k).
\end{split}\end{equation*}
In either case, we have $\big(E_{i,j}\cdot u\big)(le_{\bar i}+e_{i_1}+\dots+e_{i_k}-ke_{\bar i})=0$. Thus, $E_{i,j}\cdot u=0$.

\textbf{Claim 2.} $E_{m,m+1}\cdot u=0$.

For any $i'_1,\dots,i'_{k'}\in\{\bar i+1,\dots,m\}$,
\begin{equation*}\begin{split}
&E_{m,m+1}\cdot Y_{le_{\bar i}+e_{i'_1}+\dots+e_{i'_{k'}}-k'e_{\bar i}}v(i'_1,\dots,i'_{k'})\\
=&(\sum_{n=1}^{k'}\delta_{i'_n,m})(1-l+\mu_m+\sum_{n=1}^{k'}\delta_{i'_n,m})Y_{le_{\bar i}+e_{i'_1}+\dots+e_{i'_{k'}}-k'e_{\bar i}-e_m}v(i'_1,\dots,i'_{k'})\\
&+\sum_{s=\bar i}^{m-1}(le_{\bar i}+e_{i'_1}+\dots+e_{i'_{k'}}-k'e_{\bar i})_sY_{le_{\bar i}+e_{i'_1}+\dots+e_{i'_{k'}}-k'e_{\bar i}-e_s}E_{m,s}\cdot v(i'_1,\dots,i'_{k'})
\end{split}\end{equation*}
by Lemma \ref{Lem4.1} (5).
Then any element in $S_{E_{m,m+1}\cdot Y_{le_{\bar i}+e_{i'_1}+\dots+e_{i'_{k'}}-k'e_{\bar i}}v(i'_1,\dots,i'_{k'})}$ is of the form $le_{\bar i}+e_{i_1}+\dots+e_{i_k}-ke_{\bar i}-e_j$, where $j\in\{\bar i,i_1,\dots,i_k\}$. Moreover, if
$$le_{\bar i}+e_{i_1}+\dots+e_{i_k}-ke_{\bar i}-e_j\in S_{E_{m,m+1}\cdot Y_{le_{\bar i}+e_{i'_1}+\dots+e_{i'_{k'}}-k'e_{\bar i}}v(i'_1,\dots,i'_{k'})},$$
then one of the following cases is true:
\begin{itemize}
\item[(1)] $e_{i'_1}+\dots+e_{i'_{k'}}-k'e_{\bar i}=e_{i_1}+\dots+e_{i_k}-e_j+e_s-ke_{\bar i},\ j\in\{i_1,\dots,i_k\},s\in\{\bar i,\dots,m\}$;
\item[(2)] $e_{i'_1}+\dots+e_{i'_{k'}}-k'e_{\bar i}=e_{i_1}+\dots+e_{i_k}+e_s-ke_{\bar i},\ j=\bar i,s\in\{\bar i,\dots,m\}$.
\end{itemize}
Suppose $j=i_r$ for some $r\in\{1,\dots,k\}$. For $s\in\{\bar i+1,\dots,m\}$ we have
\begin{equation*}\begin{split}
&\big(E_{m,m+1}\cdot Y_{le_{\bar i}+e_{i_1}+\dots+e_{i_k}-e_j+e_s-ke_{\bar i}}v(i_1,\dots,i_{r-1},i_{r+1},\dots,i_k,s)\big)\\
&(le_{\bar i}+e_{i_1}+\dots+e_{i_k}-ke_{\bar i}-e_j)\\
=&\left\{\begin{aligned}
(\sum_{n=1}^k\delta_{i_n,m}-\delta_{j,m}+1)(1-l+\mu_m+\sum_{n=1}^k\delta_{i_n,m}-\delta_{j,m}+1)&\\
\cdot v(i_1,\dots,i_{r-1},i_{r+1},\dots,i_k,m), &\text{ if } s=m; \\
(\sum_{n=1}^k\delta_{i_n,s}-\delta_{j,s}+1)(1+\sum_{n=1}^k\delta_{i_n,m}-\delta_{j,m})&\\
\cdot v(i_1,\dots,i_{r-1},i_{r+1},\dots,i_k,m), &\text{ if }\ \bar i<s<m,
\end{aligned}\right.
\end{split}\end{equation*}
from Lemma \ref{Lem3.6}. For $s=\bar i$, we have
\begin{equation*}\begin{split}
&\big(E_{m,m+1}\cdot Y_{le_{\bar i}+e_{i_1}+\dots+e_{i_k}-e_j+e_s-ke_{\bar i}}v(i_1,\dots,i_{r-1},i_{r+1},\dots,i_k)\big)\\
&(le_{\bar i}+e_{i_1}+\dots+e_{i_k}-ke_{\bar i}-e_j)\\
=&(\mu_{\bar i}-\mu_{\bar i+1}-k+1)(1+\sum_{n=1}^k\delta_{i_n,m}-\delta_{j,m})v(i_1,\dots,i_{r-1},i_{r+1},\dots,i_k,m).
\end{split}\end{equation*}
Then
\begin{equation*}\begin{split}
&\big(E_{m,m+1}\cdot u)(le_{\bar i}+e_{i_1}+\dots+e_{i_k}-ke_{\bar i}-e_j)\\
=&\big((\sum_{n=1}^k\delta_{i_n,m}-\delta_{j,m}+1)(1-l+\mu_m+\sum_{n=1}^k\delta_{i_n,m}-\delta_{j,m}+1)\\
&+\sum_{s=\bar i+1}^{m-1}(\sum_{n=1}^k\delta_{i_n,s}-\delta_{j,s}+1)(1+\sum_{n=1}^k\delta_{i_n,m}-\delta_{j,m})\\
&+(\mu_{\bar i}-\mu_{\bar i+1}-k+1)(1+\sum_{n=1}^k\delta_{i_n,m}-\delta_{j,m})\big)v(i_1,\dots,i_{r-1},i_{r+1},\dots,i_k,m)\\
=&0.
\end{split}\end{equation*}

Similarly, we can verify that $\big(E_{m,m+1}\cdot u)(le_{\bar i}+e_{i_1}+\dots+e_{i_k}-ke_{\bar i}-e_j)=0$ for $j=\bar i$. Thus, $E_{m,m+1}\cdot u=0$.
\end{proof}

\begin{theorem}\label{thmmain}
$M_\p(V)$ is not simple if and only if at least one of the following cases happens:
\begin{itemize}
\item[(1)] $\bar i\neq m,\mu_{\bar i}+m+1-\bar i\in\N$ and $\mu_{\bar i}-\mu_{\bar i+1}\in\C\setminus\{0,\dots,l-1\}$;
\item[(2)] $\mu_m+1\in\N$.
\end{itemize}
\end{theorem}
\begin{proof}
According to Lemmas \ref{Lem4.9} and \ref{Lem4.10}, if either condition (1) or (2) occurs, there exists a highest weight vector in $M_\p(V)$ that is not in $\C v_\mu$. Consequently, $M_\p(V)$ can not be simple.

Now, assuming that $M_\p(V)$ is not simple, let $w$ and $\bar\alpha$ be as previously defined. Then $\bar\alpha\neq \bf 0$. Since $\mu_1=\dots=\mu_{\bar i},\ \mu_{\bar i+1}=\dots=\mu_m$, we have $\bar\alpha_1=\dots=\bar\alpha_{\bar i-1}=\bar\alpha_{\bar i+1}=\dots=\bar\alpha_{m-1}=0$ by Corollary \ref{Cor4.5}.

By Lemma \ref{Lem4.6} (2), if $\bar i\neq m$ and $\bar\alpha_{\bar i}\neq 0$, then $\bar\alpha_{\bar i}+\bar i-\mu_{\bar i}-m-1=0$. So $l=\mu_{\bar i}+m+1-\bar i=\bar\alpha_{\bar i}\in\N$. Then by Corollary \ref{Cor4.5} we have $\mu_{\bar i}-\mu_{\bar i+1}\in\C\setminus\{0,\dots,l-1\}$.

If $\bar\alpha_{\bar i}=0$, then $\bar i\neq m$ and $\bar\alpha_m\neq 0$. Again by Lemma \ref{Lem4.6} (2), $\bar\alpha_m-\mu_m-1=0$. So $\mu_m+1=\bar\alpha_m\in\N$.
\end{proof}

\section{A restriction of tensor modules over $W_m$ to $\sl(m+1)$}
Let $d=\sum_{s=1}^mt_s\pt{s}\in W_m$ and $F(P,W)$ be a tensor module over $W_m$ as defined in Subsection \ref{ssectensorm}.
It is well known that there is an embedding from $\s$ into $W_m$ given by
$$E_{i,j}\mapsto t_i\pt{j},\ E_{m+1,i}\mapsto -\pt i,\ E_{i,m+1}\mapsto t_id,\ E_{i,i}-E_{m+1,m+1}\mapsto d+t_i\pt i$$
where $1\leq i\neq j\leq m$. Denote the image of this embedding by $s'$, which is a subalgebra of $W_m$ and is spanned by $\{t_i\pt{j},\ t_id,\ \pt i,\mid 1\leq i,j\leq m\}$. Then the tensor module $F(P,W)$ is an $s'$-module by restriction.

Fix $\lambda\in(\C^*)^m$.
Define another linear isomorphism $\sigma_\lambda:\s\rightarrow\s'$ by
$$\sigma_\lambda(E_{i,j})=\lambda_i^{-1}\lambda_jt_i\pt j-\lambda_j\pt{j},\ \sigma_\lambda(E_{m+1,i})=-\lambda_i\pt{i},\ \sigma_\lambda(E_{i,m+1})=L_i,$$
$$\sigma_\lambda(E_{i,i}-E_{m+1,m+1})=\sum_{s=1}^m(t_{s}\pt{s}-\lambda_s\pt{s})+(t_{i}\pt{i}-\lambda_i\pt{i}),$$
where $L_i=\lambda_i^{-1}t_id-d-\sum\limits_{s=1}^m(\lambda_s\lambda_i^{-1}t_i\pt s-\lambda_s\pt s)$ and $1\leq i\neq j\leq m$.
\begin{lemma}
The map $\sigma_\lambda$ is an isomorphism of Lie algebras.
\end{lemma}
\begin{proof}
This is verified directly.
\end{proof}
Define the $K_m$-module $\Omega(\lambda):=\C[x_1,\dots,x_m]$ by
$$t_i^{\pm 1}\cdot f(x_1,\dots,x_m)=\lambda_i^{\pm 1}f(x_1,\dots,x_i\mp 1,\dots,x_m),$$
$$t_i\pt{i}\cdot f(x_1,\dots,x_m)=x_if(x_1,\dots,x_m),$$
where $f(x_1,\dots,x_m)\in\Omega(\lambda)$ and $1\leq i\leq m$. It is easy to verify that $\Omega(\lambda)$ is a simple $K_m$-module. Let
$X_j^i=(x_i+1)\cdots(x_i+j)\in\Omega(\lambda)$ for $1\leq i\leq m$ and $j\in\N$, and $X_j^i=1$ for $j\in-\Z_+$. For any $\eta=(\eta_1,\ldots,\eta_m)\in \Z^m$, let $X_\eta=X_{\eta_1}^1\cdots X_{\eta_m}^m$. Then $X_\eta,\eta\in\Z_+^m$ form a basis of $\Omega(\lambda)$.
Set $P=\Omega(\lambda)$, then consider the actions of the generators in $\s'$ on $F(\Omega(\lambda), W)$.
\begin{lemma}\label{Lem10}
Let $\eta\in\Z_+^m, v\in W$ and $1\leq i\neq j\leq m$. Then
\begin{itemize}
\item[(1)] $(\lambda_i^{-1}\lambda_jt_i\pt{j}-\lambda_{j}\pt{j})\cdot(X_\eta\otimes v)=-\eta_iX_{\eta-e_i+e_j}\otimes v+\lambda_i^{-1}\lambda_j
    X_\eta\otimes E_{i,j}\cdot v$;
\item[(2)] $(\lambda_i\pt{i}-t_i\pt{i})\cdot(X_\eta\otimes v)=(\eta_i+1)X_\eta\otimes v-X_\eta\otimes E_{i,i}\cdot v$;
\item[(3)] $L_i\cdot (X_\eta\otimes v)=\eta_i(|\eta|+m)X_{\eta-e_i}\otimes v+\lambda_i^{-1}\sum_{s=1}^m(-\eta_s)\lambda_sX_{\eta-e_s}\otimes E_{i,s}\cdot v-\eta_iX_{\eta-e_i}\otimes(\sum_{s=1}^mE_{s,s})\cdot v.$
\end{itemize}
\end{lemma}
\begin{proof}
(1)
\begin{equation*}\begin{split}
(\lambda_i^{-1}\lambda_jt_i\pt{j}-\lambda_j\pt{j})\cdot(X_\eta\otimes v)=&x_i^{1-\delta_{\eta_i,0}}X_{\eta-e_i+e_j}\otimes v+\lambda_i^{-1}\lambda_jX_\eta\otimes E_{i,j}\cdot v-X_{\eta+e_j}\otimes v\\
=&-\eta_iX_{\eta-e_i+e_j}\otimes v+\lambda_i^{-1}\lambda_jX_\eta\otimes E_{i,j}\cdot v.
\end{split}\end{equation*}

(2)
$$(\lambda_i\pt{i}-t_i\pt{i})\cdot(X_\eta\otimes v)=X_{\eta+e_i}\otimes v-x_iX_\eta\otimes v-X_\eta\otimes E_{i,i}\cdot v=(\eta_i+1)X_\eta\otimes v-X_\eta\otimes E_{i,i}\cdot v.$$

(3)
\begin{eqnarray*}
&&L_i\cdot(X_\eta\otimes v)=\sum_{\substack{s=1\\s\neq i}}^mx_sx_i^{1-\delta_{\eta_i,0}}X_{\eta-e_i}\otimes v+(x_i-1)x_i^{1-\delta_{\eta_i,0}}X_{\eta-e_i}\otimes v\\
&&+\lambda_i^{-1}\sum_{s=1}^m x_s^{1-\delta_{\eta_s,0}}X_{\eta-e_s}\otimes E_{i,s}\cdot v+x_i^{1-\delta_{\eta_i,0}}X_{\eta-e_i}\otimes (\sum_{s=1}^m E_{s,s})\cdot v\\
&&-\sum_{s=1}^m(x_sX_\eta\otimes v+X_\eta\otimes E_{s,s}\cdot v)+(\eta_i+1)X_\eta\otimes v-X_\eta\otimes E_{i,i}\cdot v\\
&&-\sum_{\substack{s=1\\s\neq i}}^m(-\eta_iX_{\eta-e_i+e_{s}}\otimes v+\lambda_i^{-1}\lambda_{s}X_\eta\otimes E_{i,s}\cdot v)\\
&=&\sum_{s=1}^mx_s(-\eta_i)X_{\eta-e_i}\otimes v-x_i^{1-\delta_{\eta_i,0}}X_{\eta-e_i}\otimes v+\lambda_i^{-1}\sum_{s=1}^m(-\eta_s)X_{\eta-e_s}\otimes E_{i,s}\cdot v\\
&&-\eta_iX_{\eta-e_i}\otimes (\sum_{s=1}^mE_{s,s})\cdot v+X_\eta\otimes v+\eta_i\sum_{s=1}^mY_{\eta-e_i+e_{s}}\otimes v\\
&=&\eta_i\sum_{s=1}^m(\eta_s+1)X_{\eta-e_i}\otimes v+\lambda_i^{-1}\sum_{s=1}^m(-\eta_s)\lambda_sX_{\eta-e_s}\otimes E_{i,s}\cdot v-\eta_iX_{\eta-e_i}\otimes(\sum_{s=1}^mE_{s,s})\cdot v\\
&=&\eta_i(|\eta|+m)X_{\eta-e_i}\otimes v+\sum_{s=1}^m(-\eta_s)X_{\eta-e_s}\otimes E_{i,s}\cdot v-\eta_iX_{\eta-e_i}\otimes(\sum_{s=1}^mE_{s,s})\cdot v.
\end{eqnarray*}
\end{proof}
For a $\gl(m)$-module $W$, let $\widetilde{W}$ be the same vector space as $W$, and define the $\mathfrak l$-module $\widetilde{W}$ as follows:
$$E_{i,j}\cdot v=\lambda_i^{-1}\lambda_jE_{i,j}\cdot v,$$
$$(E_{i,i}-E_{m+1,m+1})\cdot v=(\sum_{s=1}^mE_{s,s}+E_{i,i}-m-1)\cdot v$$
for $1\leq i\neq j\leq m$, where the first ``$\cdot$" is the action on $\widetilde{W}$ and the second ``$\cdot$" is the action on $W$.

\begin{theorem}\label{isoFPM}
$M_\p(\widetilde W)\cong F(\Omega(\lambda),W)^{\sigma_\lambda}$ as $\s$-modules.
\end{theorem}
\begin{proof}
From Lemma \ref{Lem10} (1)(2), we can see
$$\sigma_\lambda(\mathfrak l)\cdot (\C\otimes W)\subset\C\otimes W\subset F(\Omega(\lambda),W).$$
It is easy to verify that the $\mathfrak l$-submodule $(\C\otimes W)^{\sigma_\lambda}$ of the $\s$-module $F(\Omega(\lambda),W)^{\sigma_\lambda}$ is isomorphic to $\widetilde W$. By Lemma \ref{Lem10} (3),
$$\sigma_\lambda(E_{i,m+1})\cdot(\C\otimes W)=L_i\cdot(\C\otimes W)=0$$
for any $1\leqslant i\leqslant m$.
Thus, there is an $\s$-module homomorphism $f:M_\p(\widetilde W)\rightarrow F(\Omega(\lambda),W)^{\sigma_\lambda}$ such that $f(v)=1\otimes v$ for any $v\in W$. Note that $\sigma_\lambda(E_{m+1,i})=-\lambda_i\pt{i}$ acts on $F(\Omega(\lambda),W)$ freely. So the homomorphism $f$ is in fact an isomorphism.
\end{proof}
By this theorem, the $\s$-module $F(\Omega(\lambda),W)$ is equivalent to the generalized Verma module $M_\p(\widetilde W)$ with respect to the algebra isomorphism $\sigma_\lambda$.

\begin{corollary}
Let $\lambda\in(\C^*)^m$ and $\mu\in\C^m$ with
$$\mu_1=\dots=\mu_{\bar i},\ \mu_{\bar i+1}=\dots=\mu_m$$
for some $\bar i\in\{1,\dots,m\}$, where $\bar i=m$ if $\mu_1=\dots=\mu_m$. Then the $\s$-module $F(\Omega(\lambda),L(\mu))$ is not simple if and only if least one of the following cases happens:
\begin{itemize}
\item[(1)] $\bar i\neq m,\mu_{\bar i}+|\mu|-\bar i\in\N$ and $\mu_{\bar i}-\mu_{\bar i+1}\in\C\setminus\{0,\dots,\mu_{\bar i}+|\mu|-\bar i-1\}$;
\item[(2)] $\mu_m+|\mu|-m\in\N$.
\end{itemize}
\end{corollary}
\begin{proof}
Let $W=L(\mu)$. By Theorem \ref{isoFPM}, the $\s$-module $F(\Omega(\lambda),W)$ is not simple if and only if $M_\p(\widetilde W)$ is not simple. By the construction of the $\mathfrak l$-module $\widetilde W$, the highest weight of $\widetilde V$ is $\mu+(|\mu|-m-1)\sum_{s=1}^me_s$. Then by Theorem \ref{thmmain}, $M_\p(\widetilde W)$ is not simple if and only if at least one of the following cases happens:
\begin{itemize}
\item[(1)] $\bar i\neq m,\mu_{\bar i}+|\mu|-\bar i\in\N$ and $\mu_{\bar i}-\mu_{\bar i+1}\in\C\setminus\{0,\dots,\mu_{\bar i}+|\mu|-\bar i-1\}$;
\item[(2)] $\mu_m+|\mu|-m\in\N$.
\end{itemize}
\end{proof}
\noindent {\bf Acknowledgments.}
\noindent The authors would like to thank the professor Rencai L\"{u} for his valuable suggestions in preparation of this paper. This work was supported by NSF of China (Grant Nos. 12301037 and 12271383). Part of this work was finished during Y. Xue's visit to Chern Institute of Mathematics.

\vspace{0.4cm} \noindent


\begin{thebibliography}{00}
\bibitem{BJ} Z. Bai, J. Jiang, Gelfand–Kirillov dimension and reducibility of Scalar Generalized Verma Modules for Classical Lie Algebras, Acta Math. Sin. (Engl. Ser.), {\bf 40(3)} (2024) 658-706.
\bibitem{BXX} Z. Bai, W. Xiao, J. Jiang, Gelfand–Kirillov dimensions and associated varieties of highest weight modules, Int. Math. Res. Not., {\bf 10} (2023), 8101-8142.

\bibitem{CF} J. Coleman, V. Futorny, Stratified $L$-modules, J. Algebra, {\bf 163(1)} (1994) 219–234.

\bibitem{CLNZ} Y. Cai, G. Liu, J. Nilsson, K. Zhao, Generalized Verma modules over $\sl(n+2)$ induced from $U(\h)$-free $\sl(n+1)$-modules, J. Algebra, {\bf 502} (2018) 146-162.

\bibitem{DMP} I. Dimitrov, O. Mathieu, I. Penkov, On the structure of weight modules, Trans. Amer. Math. Soc., {\bf 352(6)} (2000) 2857–2869.

\bibitem{Fe} S. Fernando, Lie algebra modules with finite-dimensional weight spaces, I, Trans. Amer. Math. Soc., {\bf 322(2)} (1990) 757–781.

\bibitem{FM} V. Futorny, V. Mazorchuk, Structure of $\alpha$-stratified modules for finite-dimensional Lie algebras, I, J. Algebra, {\bf 183} (1996) 456-482.

\bibitem{Fu1} V. Futorny, The weight representations of semisimple finite dimensional Lie algebras, Ph.D. thesis, Kiev University (1987).

\bibitem{Fu2} V. Futorny, A generalization of Verma modules and irreducible representations of the Lie algebra $\sl(3)$, Ukrain. Mat. Zh., {\bf 38(4)} (1986) 492–497.

\bibitem{GL} H. Garland, J. Lepowsky, Lie algebra homology and the Macdonald-Kac formulas, Invent. Math., {\bf 34(1)} (1976) 37–76.

\bibitem{GLLZ} X. Guo, G. Liu, R. L\"u, K. Zhao, Simple Witt modules that are finitely generated over the Cartan subalgebra, Moscow Math. J., {\bf 20(1)} (2020) 43-65.

\bibitem{GLZ} X. Guo, X. Liu, F. Zhang, Simple $\sl_{d+1}$-modules from Witt algebra modules, Forum Math., {\bf 36(2)} (2024) 447–467.

\bibitem{J} J. Jantzen, Kontravariante formen auf induzierten darstellungen halbeinfacher Lie-algebren, Math. Ann., {\bf 226(1)} (1977) 53-65.

\bibitem{KM1} O. Khomenko, V. Mazorchuk, Generalized verma modules induced from $\sl_2$ and associated verma modules, J. Algebra, {\bf 242} (2001) 561–576.

\bibitem{KM2} O. Khomenko, V.Mazorchuk, Structure of modules induced from simple modules with minimal annihilator, Canad. J. Math., {\bf 56(2)} (2004) 293–309.

\bibitem{La} T. Larsson, Conformal fields: A class of representations of vect(N), Internat. J. Modern Phys. A, {\bf 7(26)} (1992) 6493–6508.

\bibitem{L} J. Lepowsky, Generalized Verma modules, the Cartan-Helgason theorem, and the Harish-Chandra homomorphism, J. Algebra, {\bf 2} (1977) 470–495.

\bibitem{LLZ} G. Liu, R. L\"u, K. Zhao, Irreducible Witt modules from Weyl modules and $\gl_n$-modules, J. Algebra, {\bf 511} (2018) 164–181.

\bibitem{M} V. Mazorchuk, Generalized Verma Modules, Mathematical Studies Monograph Series, vol. 8. Lviv, VNTL Publishers (2000).

\bibitem{Mc} E. McDowell, A module induced from a Whittaker module, Proc. Amer. Math. Soc., {\bf 118} (1993) 349-354.

\bibitem{MaS} V. Mazorchuk, C. Stroppel, Categorification of (induced) cell modules and the rough structure of generalised Verma modules, Adv. Math., {\bf 219} (2008) 1363–1426.

\bibitem{MS} D. Miličić, W. Soergel, The composition series of modules induced from Whittaker modules, Comment. Math. Helv., {\bf 72(4)} (1997) 503–520.

\bibitem{R} A. Rocha-Caridi, Splitting criteria for $\g$-modules induced from a parabolic and a Bernstein-Gelfand-Gelfand resolution of a finite-dimensional, irreducible $\g$-module, Trans. Amer. Math. Soc., {\bf 262} (1980) 335-366.

\bibitem{Se} A. Semikhatov, Past the highest-weight, and what you can find there, in: New developments in quantum field theory, Zakopane, (1997) 329–339.

\bibitem{S} G. Shen, Graded modules of graded Lie algebras of Cartan type (I)-Mixed products of modules, Sci. Sinica Ser. A, {\bf 29(6)} (1986) 570–581.

\bibitem{TZ} H. Tan, K. Zhao, Irreducible modules over Witt algebras $\mathcal{W}_n$ and over $\sl_{n+1}(\C)$, Algebr. Represent. Theory, {\bf 21} (2018) 787–806.

\bibitem{V} D. Verma, Structure of certain induced representations of complex semisimple lie algebras, Ph. D. thesis, Yale University (1966).

\bibitem{ZX} Y. Zhao, X. Xu, Generalized projective representations for $\sl(n +1)$, J. Algebra, {\bf 328} (2011) 132-154.
\end{thebibliography}
\end{document}